\documentclass{amsart}
\newtheorem{teo}{Theorem}[section]

\newtheorem{prop}[teo]{Proposition}
\newtheorem{corollario}[teo]{Corollary}
\theoremstyle{definition} 
\newtheorem{definiz}[teo]{Definition} 
\newtheorem{example}[teo]{Example}
\newenvironment{proof1}{\noindent\textit{Proof}}{\hfill$\square$\medskip\\ }

\usepackage{epsfig,amsfonts,amssymb}

\theoremstyle{remark} 
\newtheorem{nota}[teo]{Remark}
\numberwithin{equation}{section}

\begin{document}
\title[On Closed Invariant Sets]{On Closed Invariant Sets in Local Dynamics}
 \author[C. Bisi]
{Cinzia Bisi $^\ast$}
\date{August 13.th 2008}
\thanks{\rm $^\ast$ Partially supported by Progetto MURST di
Rilevante Interesse Nazionale {\it Propriet{\`a} geometriche
delle variet{\`a} reali e complesse}, by GNSAGA and by INDAM}
\address{Dipartimento di Matematica, Universit\'a della Calabria, Ponte Bucci, Cubo 30b, 87036, Arcavacata di Rende (CS), Italy}
\vspace{1cm}
\email{bisi@math.unifi.it}
\email{bisi@mat.unical.it}

\begin{abstract}
We investigate the dynamical behaviour of a holomorphic map on a $f-$invariant subset $\mathcal{C}$ of $U,$ where
$f:U \to \mathbb{C}^k.$ We study two cases: when $U$ is 
an open, connected and polynomially convex subset of $\mathbb{C}^k$ and $\mathcal{C} \subset \subset U,$ closed in $U,$ and when $\partial U$
has a p.s.h. barrier at each of its points and $\mathcal{C}$ is not relatively compact in $U.$ In the second part of the paper,  we prove a Birkhoff's type Theorem for holomorphic maps in several complex variables, i.e. given an injective holomorphic map $f,$ defined in a neighborhood of $\overline{U},$ with $U$ star-shaped and $f(U)$ a Runge domain, we prove the existence of a unique, forward invariant, maximal, compact and connected subset of $\overline{U}$ which touches $\partial U.$
  
%In this paper we give a new version of a recent result of Fornaess and Stensones,\cite{FS}, %about the non existence of dense orbits for a germ of holomorphic self-map of $\mathbb{C}^k.$
%We prove that, given an holomorphic map $f:U \to \mathbb{C}^k,$ with $U$ an open, polynomially %convex subset of $\mathbb{C}^k,$ $f$ doesn't have a dense orbit. Therefore the %Fornaess-Stensones' condition that $U$ has a Lipschitz boundary is not necessary, as pointed %out from the same authors in their paper, \cite{FS}.
\end{abstract}

\maketitle

\section{Introduction}

Let $f:U \to \mathbb{C}^k$ be a holomorphic map. Here $U$ is an open, connected and bounded (or hyperbolic) subset in $\mathbb{C}^k.$ Since the semi-local holomorphic dynamics isn't well understood yet, specially when $k>2,$
we describe the dynamical behaviour of $f$ on a $f-$invariant subset $\mathcal{C}$ of $U$ in two different cases: \\
\begin{itemize}
\item[a)] when $\mathcal{C} \subset \subset U,$ closed in $U,$ and $U$ is polynomially convex; \\
\item[b)] when $\mathcal{C}$ isn't relatively compact in $U$ and every point in $\partial U$ has a p.s.h. barrier.
\end{itemize}
When there is a recurrent component $V$ in the interior of the polynomially convex hull of $\mathcal{C}$ in case $a)$ or 
in the interior of $\overline{\mathcal{C}}$ in case $b),$ we prove that the dynamical behaviour on $V$ is of three types:
\begin{itemize}
\item[1.] $V$ is the basin of attraction of an attractive periodic orbit; \\
\item[2.] $V$ is a Siegel domain;\\
\item[3.] if $h$ is a limit of a subsequence of $\{f^n\}_{n \in \mathbb{N}},$ then $0 < \textrm{rank} \,(h) < k.$ 
\end{itemize}
In particular when $\mathcal{C}$ is a closed orbit or a countable union of closed orbits, we prove that $\mathcal{C}$ cannot have a non-empty interior with a recurrent point. This has been proved by Fornaess-Stensones in \cite{FS} when $U$ has a Lipschitz boundary; here it is proved in a different situation, i.e. when $U$ is polynomially convex or with a p.s.h. barrier at each boundary point, then $U$ has not necessarily Lipschitz boundary. \\
In the second part of the paper, see section \ref{B}, we give a version of Birkhoff's Theorem which was originally stated for surface transformations $f$ having a Lyapunov unstable fixed point $p$ for $f$ or for $f^{-1}.$ Under these hypotheses Birkhoff has shown, \cite{Bir}, the existence, in each neighborhood $U$ of $p,$ of a compact set $K_{+}$ (or $K_{-}$) which is positive (or negative) invariant by $f$ and touching the boundary of $U.$ In this general setting there is no forward and backward invariant compact set with this property.
\\In the same spirit, our theorem \ref{birkhoff} asserts that if $f:U \to \mathbb{C}^k$ is a holomorphic injective map of $\mathbb{C}^k$ such that $f(0)=0,$ with $U$ bounded and star-shaped and $f(U)$ a Runge domain, then there exists a unique, maximal, compact, connected set $K$ such that:
\begin{itemize}
\item[1.] $0 \in K \subset \overline{U};$ \\
\item[2.] $K \cap \partial U \neq \emptyset;$ \\
\item[3.] $f(K) \subset K.$
\end{itemize}
In general, this compact set $K$ isn't totally invariant: 
we will give an example, see example \ref{exhen}. So the several variables analogue of R. Perez-Marco's {\it hedgehogs}, \cite{PM}, does not hold: in the one variable case the compact is totally invariant and touches the boundary, \cite{PM}.

\section{Preliminaries}
We recall some definitions and fix our notations. \\
Let $K$ be a compact set of $\mathbb{C}^k,$ then the polynomially convex hull of $K$ is defined as: \\
$\hat{K}_{\mathcal{P}}= \{ z \in \mathbb{C}^k \,\,\,\, | \,\,\,\, |p(z)| \le \sup\limits_{\zeta \in K} |p(\zeta)| \,\,\,\,\,
\forall \,\, p \,\,\, \mathrm{polynomial} \}.$ \\
A compact set $K$ is {\it polynomially convex} if $K=\hat{K}_{\mathcal{P}}.$
\begin{definiz}
An open set $U$ in $\mathbb{C}^k$ is {\it polynomially convex} if, for every compact $K$ in $U,$ $\hat{K}_{\mathcal{P}} \subset \subset U.$ 
\end{definiz}
For example, the geometrically convex open sets of $\mathbb{C}^k$ are polynomially convex in $\mathbb{C}^k.$ 
The property of being polynomially convex is not invariant by biholomorphisms, as Wermer showed, see Gunning's book, \cite{G}, page 46.\\
If $K$ is polynomially convex, each holomorphic function on a neighborhood of $K$ is the uniform limit on $K$ of polynomials; in the same way if $\rho$ is p.s.h. and continuous on $U,$ polynomially convex open set, then it is the uniform limit on the compact sets of $U$ of p.s.h. functions of $\mathbb{C}^k.$ \\
A consequence, when $U$ is polynomially convex, is that convexity with respect to p.s.h. functions in $U$ is the same as polynomial convexity. \\
If $K$ is polynomially convex and compact in $U,$ there exists $\rho_1$ p.s.h. and continuous on $\mathbb{C}^k,$ $K=\{ \rho_1 \le 0 \}$ and $\rho_1 \ge 1$ on a neighborhood of $\mathbb{C}^k \setminus U.$ \\
\begin{definiz}
A domain $U$ is {\it Runge} if each holomorphic function on $U$ can be approximated by polynomials, uniformly on compact subsets of $U.$
\end{definiz}
In particular any polynomially convex open set is a Runge domain, \cite{G}.\\
It is possible to construct Runge domains such that the interior of $\overline{U}$ is not equal to $U:$ for example $U = \{ w \in \mathbb{C}^k : |w|< \exp(-\varphi) \}$ with $\varphi$ subharmonic on the unit disc, $\varphi =0$ on a dense set of $\Delta$ and $\varphi \ge 0$ and non identically zero, in particular it doesn't have Lipschitz boundary.

\section{Invariant Sets}

\subsection{$f-$Invariant Relatively Compact Subsets.}
\vspace{0.20cm} 
Let $f:U \to \mathbb{C}^k$ be a holomorphic map with $U \subset \subset \mathbb{C}^k$ or $U$ Kobayashi hyperbolic. We assume that $U$ is an open, connected and polynomially convex set.
We say that a closed set $\mathcal{C}$ is $f-$invariant if $f(\mathcal{C}) \subset \mathcal{C}.$
\begin{prop} \label{prop2}
Let $\mathcal{C} \subset  \subset U$ be a closed $f-$invariant set, then $\hat{\mathcal{C}}_{\mathcal{P}}$ is $f-$invariant.
\end{prop}
\begin{proof1}\\
By hypothesis, $\mathcal{C} \subset \subset U.$
Choose $z_0 \in \hat{\mathcal{C}}_{\mathcal{P}}$ and suppose $f(z_0) \notin
\hat{\mathcal{C}}_{\mathcal{P}}.$ Then there is a p.s.h. smooth function $\rho_0$ in $\mathbb{C}^k,$ such that $\rho_0 \le 0$ on $\hat{\mathcal{C}}_{\mathcal{P}}$ and $\rho_0(f(z_0)) >1.$ \\
The function $\rho_0 \circ f$ is p.s.h. on $U,$ $\rho_0 \circ f \le 0$ on $\mathcal{C}$ and $\rho_0 \circ f$ is also p.s.h. on the holomorphic hull of $\mathcal{C}$ with respect to $U,$ which is the same as $\hat{\mathcal{C}}_{\mathcal{P}}.$ It follows, by Maximum Principle, that $\rho_0 (f(z_0)) \le 0,$ which is  a contradiction. 
\end{proof1}
\begin{definiz}
A connected component $\Omega \subset U,$ of the set of points where $\{ f^n \}_{n \in \mathbb{N}}$ is equicontinuous, is {\it recurrent} if there exists $p_0 \in \Omega$ such that $f^{n_i} (p_0)$ is relatively compact in $\Omega$ for some subsequence $n_i,$ i.e. if $\Omega$ contains a recurrent point $p_0.$
\end{definiz}
\begin{prop} \label{prop3}
If $V=Int(\hat{\mathcal{C}}_{\mathcal{P}}) \neq \emptyset$ then the sequence 
$\{f^n\}_{n \in \mathbb{N}}$ defined on $V$ is a normal family and if $V$ has a recurrent component $U$ then there are three possibilities:
\begin{itemize}
\item[i)] $f$ has an attractive periodic orbit,\\
\end{itemize}
\begin{itemize}
\item[ii)] there is a Siegel domain, i.e. there is $U,$ a component of $V$ and a subsequence $n_i,$ s.t. $f^{n_i}_{|U} \to Id,$
\end{itemize}
\begin{itemize}
\item[iii)] if $h$ is a limit of a subsequence of $\{f^n\}_{n \in \mathbb{N}},$ then $0 < \textrm{rank} \,(h) < k.$
\end{itemize}
\end{prop}
\begin{proof1}\\
We assume that for some $p_0,$ $f^{n_i} (p_0) \to p \in U,$ and $f^{n_i}$ converges uniformly on compact sets. We now write
$f^{n_{i+1} - n_i} \circ f^{n_i} = f^{n_{i+1}}.$ Extracting a subsequence we get a limit $h$ of $f^{n_{i+1} - n_i}$ such that
$h(p)=p,$ \cite{FSib1}. If $h$ is of rank $0,$ we show that $p$ is an attractive fixed point, \cite{FSib1}. If $h$ is of maximal rank, then we get a Siegel domain, \cite{FSib1}. The theorem of Carath\'eodory-Cartan-Kaup-Wu, see \cite{K} page 438 and \cite{N} page 66, describes the permitted eigenvalues. Otherwise for all possible $h,$ $0 < \textrm{rank} \,(h) < k.$
\\ In \cite{FSib1}, Fornaess and Sibony prove a more precise result when $f$ is an endomorphism of $\mathbb{P}^2.$ Their stronger result is valid only in dimension two.   
\hfill
\end{proof1}

\subsection{$f-$Invariant Non-Relatively Compact Subsets.}
\begin{teo} \label{rec}
Let $f:U \to \mathbb{C}^k$ be a holomorphic open map defined on $U,$ a bounded (or hyperbolic) open and connected subset of $\mathbb{C}^k.$ Assume that every point in $\partial U$
has a p.s.h. barrier, i.e. if $q \in \partial U,$ there exists a p.s.h. function $\rho_q,$ $\rho_q <0$ on $U,$ continuous such that $\lim\limits_{p \to q} \rho_q (p)=0.$ Suppose $\mathcal{C}$ is a $f-$invariant set in $U.$ 
Let $V$ be the non-empty interior of $\overline{\mathcal{C}},$ where the adherence is with respect to $U.$ We also assume that a connected component of $V,$ $W,$ contains a recurrent point $p_0.$ Then there are three possibilities for $W:$
\begin{itemize}
\item[1)] it is the basin of attraction of an attracting periodic orbit;
\end{itemize}
\begin{itemize}
\item[2)] it is a Siegel domain;
\end{itemize}
\begin{itemize}
\item[3)] if $h$ is a limit of a subsequence of $\{f^n\}_{n \in \mathbb{N}},$ on $W,$ then $0 < \textrm{rank} \, (h) < k.$
\end{itemize}
\end{teo}
\begin{proof1}\\
We start proving that the sequence $\{ f^n \}_{n \in \mathbb{N}}$ is well defined on $V.$
Since $V \subset U$ is invariant, by continuity $f(V) \subset \overline{U}:$
indeed if $p \in V$ there exists a sequence of points $p_n \in \mathcal{C}$ such that $p_n \to p$
and hence $f(p_n) \to f(p)=q \in \overline{U}.$ We show now that $f(V) \subset U.$ Suppose $q \in \partial U.$ Consider the barrier $\rho_q$ at $q.$ The function $\rho_q \circ f$ is p.s.h. and continuous on $V,$ and $\rho_q \circ f \le 0$ on $V.$ But $(\rho_q \circ f)(p)=\lim\limits_{n \to + \infty} (\rho_q \circ f)(p_n)=\lim\limits_{n \to + \infty} \rho_q (f(p_n))=0.$ Hence, by Maximum Principle, $\rho_q \circ f \equiv 0,$ i.e. $f(V) \subset (\rho_q =0) \subset \partial U.$ This is impossible because $f$ is open. 
Hence $f(V) \subset U$ and $f^n (V) \subset U,$ therefore the sequence $\{ f^n \}_{n \in \mathbb{N}}$ is normal, since $U$ is bounded.\\
Now suppose that there exists a recurrent point $p_0$ in $W,$ a connected component of $V.$ This means that there exists a sequence of $n_i \to + \infty$ s.t. $f^{n_i} (p_0) \to p_0 \in W.$ We can always suppose that
$n_{i+1} - n_i \to +\infty.$ Taking a subsequence $\{i=i(j)\}$ we can suppose that the sequence $\{ f^{n_{i+1} - n_i} \}_i$ converges uniformly on compact sets of $W$ to a holomorphic map $h:W \to \overline{U}$ s.t. $h(p_0)=p_0.$ Indeed let $p_i =f^{n_i} (p_0).$ Then $f^{n_{i+1} -n_i} (p_i)=f^{n_{i+1}} (p_0) = p_{i+1}.$ Hence $f^{n_{i+1}-n_i} (p_0)=p_{i+1} + O(|p_i -p_0|)$ so converges to $p_0$ and therefore, necessarily, $h(p_0)=p_0,$ \cite{FSib1}.
\\
% Then from an analysis similar to the one in \cite{FSib1}, $V$ is a subset either 
%of a basin of attraction of some periodic orbit, or of a Siegel domain.  
Consider all maps $h$ obtained in this way. If some $h$ is of rank $0,$ then some iterate of $f$ has $p_0$ as an attractive fixed point and $f$ has $p_0$
as an attractive periodic point.\\
If some $h$ is of maximal rank $k,$ then $W$ is a Siegel domain, otherwise all the limit maps have lower rank $r,$ $0<r<k.$
In \cite{FSib1} the authors analyze the case of holomorphic endomorphisms of $\mathbb{P}^2$ and  thanks to the restriction to the dimension $2$ and to the endomorphism case, the result there is much more precise: for example in case $iii),$ $h(W)$ is always independent of $h$ and attracts all orbits. 
\end{proof1}
\begin{nota}
If $f$ is not open it is enough to assume that $(\rho_q=0)$ does not contain the image of $f.$
\end{nota}
\begin{corollario}
Under the hypotheses of Theorem \ref{rec}, if $\overline{\mathcal{C}}$ is an invariant closed set with a dense orbit in it or a countable union
of closed invariant sets each one with a dense orbit, then the interior $V$ of $\overline{\mathcal{C}}$ doesn't contain recurrent points. 
\end{corollario}
\begin{proof1}\\
Indeed in the possible dynamical behaviours described in Theorem \ref{rec}, when  
$\overline{\mathcal{C}}$ is closed with a dense orbit cannot have interior; when we consider a countable union of closed sets with empty interior then, by Baire's theorem, the union of them is still with empty interior.  
\end{proof1}

\section{Forward invariant compact sets} \label{B}

\begin{teo}\label{birkhoff}
Let $U$ be a bounded star-shaped domain with respect to $0$ in $\mathbb{C}^k$ and let $U'$ be an open neighborhood of $\overline{U}.$
Let $f:U' \to \mathbb{C}^k,$ be a holomorphic map, $f(0)=0,$ $f$ injective on $U$ (i.e. $f:U \to f(U)$ is a biholomorphic map) and $f(U)$ is a Runge domain. Assume $f(z)= Az +O(z^2),$ with $A$ a linear invertible map and with all the eigenvalues $\lambda_j,$ for $1 \le j \le k,$ of modulus 1. 
%, since it has a determinat different from $0$ because each complex matrix is %triangularizable).
Then there exists a unique maximal connected compact set $K,$ with $0 \in K \subset \overline{U}$ s.t.
$(K \cap \partial U) \neq \emptyset$ and $f(K)\subset K.$ Futhermore $f$ is linearizable iff $0 \in Int (K).$ 
\end{teo}
\begin{proof1}\\
Define $f_{\mu_n} (z)=f(\mu_n \cdot z)$ with $\mu_n \in \mathbb{R},$ $0< \mu_n <1$ and $\mu_n \to 1$ for $n \to +\infty.$ Then
$f_{\mu_n} \to f$ uniformly on $\overline{U}$ and $|Jac(f_{\mu_n})(0)|=|\mu_n|\cdot |Jac(A)|<1$ because $|\mu_n|<1$ and $|Jac(A)|=1;$
indeed
$f_{\mu_n} (z)=\mu_n \cdot A \cdot z + O((\mu_n \cdot z)^2).$ \\
For simplicity, we call $\mu := \mu_n.$\\
Let $f_{\mu}: \dfrac{1}{\mu} \cdot U \to f(U)$ indeed $f_\mu(\dfrac{1}{\mu} \cdot U) \equiv f(U).$ Hence $f_\mu$ is a biholomorphism from a star-shaped domain $\dfrac{1}{\mu} \cdot U$ to a Runge domain $f(U)=f_{\mu}(\dfrac{1}{\mu} \cdot U).$
Now applying a result of Andersen-Lempert, \cite{AL} Theorem 2.1, to the biholomophism $f_{\mu} : \dfrac{1}{\mu} \cdot U \to f(U),$ we find a sequence of automorphisms $g_m$ of $\mathbb{C}^k,$ such that
$g_m \to f_{\mu}$ for $m \to +\infty$ uniformly on compact subsets of $\overline{U},$ i.e. the $g_m$'s converge to $f_{\mu},$ uniformly on compact sets and $g_m(0)=0$ for all $m.$ 
%as all non-affine shears
\\
Since $|Jac(f_{\mu})(0)|<1,$ then $|Jac(g_m)(0)|<1.$  \\
Hence $g_m \in Aut(\mathbb{C}^k)$ and $g_m:U \to g_m(U)$ with $0 \in U \cap g_m(U).$ \\
Let $B$ be a domain which is a homothetic of $U,$ i.e. $B=\epsilon U,$ sufficiently small s.t.
$g_m^{-1}(B) \subset U$ i.e. $0 \in B \subset (U \cap g_m(U)).$ Since the basin of attraction of $0$ for $g_m$ (i.e. $\bigcup\limits_{n \in \mathbb{N}} g_m^{-n}(B)$) is biholomorphic to $\mathbb{C}^k,$ \cite{RR}, and in particular is unbounded, there exists $n_0 \in \mathbb{N}$ s.t. $g_m^{-n_0}(B) \subset U$ but $g_m^{-(n_0+1)}(B) \not\subset U$ ($n_0 \ge 1$). \\
We consider the one-parameter family $\{ B_t \}_{t \ge 1}$ where $B_t = t \cdot B,$ \cite{PM}. 
Then we consider the $t$'s for which:
$$g_m^{-n_0} (B_t) \subset U.$$
The set is not empty because for $t=1$ the inclusion is true.
By continuity, there exists $\overline{t}$ s.t.
$$g_m^{-n_0}(B_{\overline{t}}) \subset U$$
and
$$g_m^{-n_0}(\overline{B_{\overline{t}}}) \cap (\partial U) \neq \emptyset.$$
We call $F_m : = \overline{g_m^{-n_0}(B_{\overline{t}})}.$ \\
Then $(F_m)_{m \in \mathbb{N}}$ is a sequence of compact sets in $\overline{U}$ s.t. $g_m(F_m)\subset F_m$ because \\
$g_m^{-n_0+1}(B_{\overline{t}}) \subset g_m^{-n_0}(B_{\overline{t}}):$ this follows from the description of the basin of attraction of $0.$ \\
Each $F_m$ is a connected set because it is the closure of the pre-image by a biholomorphism of a connected set.\\
By compactness of the space $\mathcal{K}_c (\overline{U})= \{ {\textrm connected} \,\,\, {\textrm compact} \,\,\, {\textrm subsets} \,\,\, {\textrm of} \,\,\, \overline{U} \},$ there exists a subsequence $(m_k)_{k \in \mathbb{N}}$ t.c.
$F_{m_k} \to K_{\mu} \in \mathcal{K}_c (\overline{U}).$ Finally we prove that $f_{\mu}(K_{\mu}) \subset K_{\mu}.$\\
We use that: 
\begin{itemize}
\item[(i)]$g_m \to f_{\mu}$ uniformly on compact subsets of $\overline{U};$ \\
\item[(ii)] $\lim\limits_{k \to +\infty} F_{m_k}=K_{\mu}.$        
\end{itemize}
Let $x \in K_{\mu},$ then we want to prove that $f_{\mu}(x) \in K_{\mu}.$ \\
Since $x \in K_{\mu},$ there exists a sequence $x_k \to x$ with $x_k \in F_{m_k}$ by $(ii).$ \\
Then $g_{m_k}(x_k) \in F_{m_k}$ and we can assume $g_{m_k}(x_k) \to y \in K_{\mu},$ by $(ii).$ \\
But $g_{m_k} \to f_{\mu}$ for $k \to +\infty$ by $(i),$ so $f_{\mu}(x)=\lim\limits_{k \to \infty} g_{m_k}(x_k)=y \in K_{\mu}.$ \\
Hence $K_{\mu}$ is $f_{\mu}-$invariant.\\
Therefore for each $\mu$ we have found a forward invariant connected compact set for $f_{\mu}$ and $K_{\mu}$ intersects $\partial U.$ Now, with an argument similar to the one already used for $\{g_m \}_{m \in \mathbb{N}}$ and $\{ F_{m_k} \}_{k \in \mathbb{N}},$ we prove that, up to considering a subsequence, $K_{\mu_n} \to K$ in the Hausdorff metric. Since $f_{\mu_n} \to f$ uniformly on compact sets, we have that $f(K) \subset K$ and $K$ touches $\partial U.$
In order to have the unique, maximal, connected, invariant compact set, it is enough to take the closure of the union of all such compact sets $K.$ Obviously, the closure of a union of $f-$invariant sets is still $f-$invariant and it is also connected
because each compact set contains $0.$
Since $K_{\mu_n}$ intersects $\partial U$ for all $\mu_n,$ also its limit $K$ in the Hausdorff topology does.
Suppose $0 \in Int (K),$ we show that $f$ is linearizable. The family $(f^n)_{n \in \mathbb{N}}$ is locally equicontinuous on $Int(K)$ and $f(0)=0.$
Following a standard trick, we define
$$
h(z):= \lim\limits_{n_j \to + \infty} \frac{1}{n_j} \sum\limits_{j=0}^{n_j -1} A^{-j} f^j(z).
$$
The limit exists in a neighborhood of zero. Indeed there is a $c>1$ such that $f^n( \mathbb{B}(0,r)) \subset \mathbb{B}(0,cr) \subset K$ for all $n.$
Then we can consider a limit map $h$ for an appropriate subsequence $n_j.$ We have $h(0)=0,$ $Jac(h)(0)=Id$ and we easily check that $h(f)=Ah.$
\end{proof1}
\vspace{0.5cm}
\begin{nota}
If we take a sequence $\mu_n >1,$ $\mu_n \to 1,$ we can prove that there exists a maximal connected compact set invariant for $f^{-1}.$ In general the forward and backward invariant compact subsets are different, as the case of H\'enon maps shows, see
Example \ref{exhen} below.
\end{nota}
\begin{nota}
We want to point out that $K$ is not necessarily a proper subset of $\overline{U},$ indeed 
if $f$ is an automorphism of the ball $\mathbb{B}^k \subset \mathbb{C}^k$ fixing $0,$ then $K = \overline{\mathbb{B}^k}.$
\end{nota}
\begin{nota}
Suppose that $f,g$ are two commuting maps satisfying all the hypotheses of Theorem \ref{birkhoff}, then they share the same maximal, compact, connected, invariant set $K \ni 0.$ \\
Indeed let $K_f$ and $K_g$ be the maximal, compact, connected invariant sets containing $0,$ for $f$ and $g$ respectively, which exist by Theorem \ref{birkhoff}. Then consider $f \circ g (K_f)= g \circ f (K_f) \subset g(K_f),$ hence $g(K_f) \subset K_f$ which implies that $K_f \subset K_g.$
Analogously, considering $g \circ f (K_g)= f \circ g (K_g),$ we can prove that $K_g \subset K_f.$
\end{nota}

\section{Examples}

In this section we are going to prove that our theorem \ref{birkhoff} is optimal, we mean that
there exist a map $f:\mathbb{B} \to \mathbb{C}^k$ which satisfy all the hypotheses of the Theorem \ref{birkhoff} such that it has a forward invariant compact and connected set containing $0$ which touches the boundary of $\mathbb{B}$ but it doesn't admit a totally invariant compact and connected  set containing $0$ which touches the boundary of $\mathbb{B}.$
\begin{example} \label{exhen}
Let $f$ be the following H\'{e}non map: \\
$$f(z,w)=(z^2+w,z).$$
Then $f(0,0)=(0,0)$
and $$Jac(f)= \left( \begin{array}{ll}
                          2z & 1 \\
                          1 & 0
                      \end{array}
                      \right) $$
                     
So, at $0,$ $\lambda_1=1,$ $\lambda_2=-1,$ i.e. $|\lambda_j|=1$ for $j=1,2.$
Clearly $f \in Aut(\mathbb{C}^2).$
From the well known study of the dynamics of $f,$ there exist the following closed invariant subsets of $\mathbb{C}^2:$ 
$$K^+_f=\{ z \in \mathbb{C}^2 \,\,\,\, | \,\,\,\,  f^n(z) \,\,\,\,{\textrm is} \,\,\, {\textrm bounded} \}$$
$$K^-_f=\{ z \in \mathbb{C}^2 \,\,\,\, | \,\,\,\,  f^{-n}(z) \,\,\,\,{\textrm is} \,\,\, {\textrm bounded} \}$$
and the following compact set of $\mathbb{C}^2$ containing $0:$
$$K=K^+_f \cap K^-_f.$$ 
Consider a ball $\mathbb{B}(0,R) \subset \mathbb{C}^2$ with $R>>1$ such that $\mathbb{B}(0,R) \supset\supset K.$ If we consider the restriction $f: \mathbb{B} (0,R) \to \mathbb{C}^2,$ by Theorem \ref{birkhoff} there exists a connected compact subset $X$ of $\mathbb{B}(0,R)$ which touches the $\partial \mathbb{B}(0,R),$ which is $f-$invariant and which contains $0.$ For any such $X,$ we have
$X \subset K^+_f,$ \cite{Sib}, because if $z \in X,$ $f^n (z)$ is bounded since $X$ is $f-$invariant and compact. Hence $X \subset (K^+_f \cap \overline{\mathbb{B} (0,R)}).$ 
It is well known from the study of the dynamics of H\'{e}non maps that:
$$dist(f^n(X),K) \to 0$$ uniformly on compact sets.
Hence there exists $n_0 \in \mathbb{N}$ such that
$dist(f^{n_0}(X),K)< \dfrac{1}{2} \cdot dist(K, \partial \mathbb{B}(0,R)).$ So $X$ cannot be at the same time forward and backward invariant i.e. $f(X) \subset X,$ but $f(X) \neq X.$ \\
If $f^{n_0}(X)$ is distant from $K$ less than $dist(\partial \mathbb{B}(0,R), K),$
then it means that $f^{n_0}(X) \subset X$ and they are different. \\
Hence, if we consider $g:=f^{n_0},$ then $g(X) \subset \subset X.$
\end{example}
\begin{example} 
In some cases it is possible that the forward and the backward invariant compact sets coincide.
For example, if in the previous example we consider a ball $\overline{\mathbb{B}}(0,r)$ which contains $K=K^+_f \cap K^-_f$ and such that
$K \cap \partial \mathbb{B}(0,r) \neq \emptyset,$ then the restriction of the H\'{e}non map
$f$ to $\mathbb{B} (0,r)$ admits a forward and backward invariant compact set $K$ which touches the boundary of $\mathbb{B} (0,r).$ 
\end{example}
\begin{nota}
Let $K$ be one of the $f-$invariant, connected and compact set of Theorem \ref{birkhoff}, and let $X= \bigcap\limits_{n \in \mathbb{N}} f^n(K).$ The set $X$ is connected because it is a decreasing intersection of connected sets, $X \ni 0,$ $X$ is compact and $f(X)=X.$ For example if $f$ is an H\'{e}non map, $X= K^+_f \cap K^-_f.$ 
\end{nota}

\bibliographystyle{amsplain}

\begin{thebibliography}{Co-MAc}
\bibitem{Ab} M. Abate {\em The residual index and the dynamics of holomorphic maps tangent to the identity}. Duke Math. J., 107, n.1 (2001), 173-207. 
\bibitem{AL} E. Andersen, L. Lempert {\em On the group of holomorphic automorphisms of $\mathbb{C}^n$}. Invent. Math. 110, 371-388 (1992).
\bibitem{Bir} G.D. Birkhoff {\em Surface transformations and their dynamical applications}.
Acta Math., 43 (1922), 1-119; Also in: Collected Mathematical Papers, Vol. II, pp.111-229.
\bibitem{C} D. Cerveau {\em Personal Communication}. 
\bibitem{DS} T.C. Dinh, N. Sibony {\em Dynamique des applications d' allure polynomiale}. J. Math. Pures Appl. (9), 82 (2003), n.4, 367-423.
\bibitem{FS} J.E. Fornaess, B. Stensones {\em Density of orbits in complex dynamics}.
Ergod. Th. and  Dynam. Sys. (2006), 26, 169-178. 
\bibitem{FSib1} J.E. Fornaess, N. Sibony {\em Classification of recurrent domains for some holomorphic maps}. Math. Ann. 301,813-820 (1995).
\bibitem{FSib2} J.E. Fornaess, N. Sibony {\em The closing lemma for holomorphic maps}. Ergod. Th. and Dynam. Sys. (1997), 17, 821-837.  
\bibitem{FR1} F. Forstneric, J.P. Rosay {\em Approximation of biholomorphic mappings by automorphisms of $\mathbb{C}^n$}. Invent. Math., 112, 323-349 (1993).
\bibitem{FR2} F. Forstneric, J.P. Rosay {\em Erratum. Approximation of biholomorphic mappings by automorphisms of $\mathbb{C}^n$}. Invent. Math., 118, 573-574 (1994).
\bibitem{G} R.C. Gunning {\em Introduction to Holomorphic Functions of Several Variables}. Vol.I: function theory. Wadsworth and Brooks/Cole Publishing Company, 1990.
\bibitem{Hak} M. Hakim {\em Analytic transformations of ($\mathbb{C}^p,0$) tangent to identity}.
Duke Math. J. 92 (1998), n.2, 403-428.
\bibitem{H} L. Hormander {\em An Introduction to Complex Analysis in Several Variables}. Third Edition, North-Holland Mathematical Library.
\bibitem{N} R. Narasimhan {\em Several Complex Variables} Chicago Lectures in Mathematics Series, 1971.
\bibitem{PM} R. Perez-Marco {\em Fixed Points and Circle Maps}. Acta Math., 179 (1997), 243-294.
\bibitem{RR} J.P. Rosay, W. Rudin {\em Holomorphic maps from $\mathbb{C}^n$ to $\mathbb{C}^n$}.
Trans. Amer. Math. Soc., 310, (1988), 47-86. 
\bibitem{Sib} N. Sibony {\em Dynamique et G\'{e}om\'{e}trie Complexes}. Panoramas et Synth\'{e}ses, n. 8, 1999, Soci\'{e}t\'{e} Math\'{e}matique de France.
\bibitem{K} S. Krantz {\em Function Theory of Several Complex Variables}. Second Edition, Wadsworth and Brooks/Cole.
\end{thebibliography}

\end{document}